\documentclass{amsart}

\usepackage{amsmath,amssymb,dsfont}

\usepackage{graphicx}  

\numberwithin{equation}{section}

\newcommand{\abs}[1]{| #1|}

\newcommand{\C}{\mathbb{C}}
\newcommand{\ddt}{\frac{d}{dt}}
\newcommand{\ep}{\epsilon}

\newcommand{\pp}{\mathbb{P}}
\newcommand{\R}{\mathbb{R}}
\newcommand{\T}{\mathbb{T}}
\newcommand{\Z}{\mathbb{Z}}
\renewcommand{\to}{\rightarrow}

\renewcommand{\leq}{\leqslant}
\renewcommand{\geq}{\geqslant}

\newcommand{\eps}{\varepsilon}

\newtheorem{thm}{Theorem}[section]

\newtheorem*{defin*}{Definition}

\newtheorem{corollary}{Corollary}[section]

\numberwithin{equation}{section}

\begin{document}

\title{Solitons and Gibbs measures for nonlinear Schr\"odinger equations}
\author{Kay Kirkpatrick}
\address{University of Illinois at Urbana-Champaign, Department of Mathematics, 1409 W. Green St., Urbana, IL 61801}
\email{kkirkpat@illinois.edu}
\thanks{Supported in part by NSF grants OISE-0730136 and DMS-1106770.}

\maketitle

\begin{abstract}
 We review some recent results concerning Gibbs measures for nonlinear Schr\"odinger equations (NLS), with implications for the theory of the NLS, including stability and typicality of solitary wave structures. In particular, we discuss the Gibbs measures of the discrete NLS in three dimensions, where there is a striking phase transition to soliton-like behavior.
\end{abstract}

\section{Introduction}

The cubic nonlinear Schr\"odinger equation (NLS) arises naturally in physics and engineering, for instance, as a model of nonlinear optics, Langmuir waves in plasmas, and Bose-Einstein condensation (BEC)~\cite{LRS, Z, BKA, FKM, W, rumpf04}. It turns out that the NLS can be rigorously derived from the physics of interacting particles in BEC via a scaling limit and statistical mechanics~\cite{ESY2, ESY,KSS}. It also turns out that the probabilistic tools of statistical mechanics are key to understanding the dynamics of nonlinear Schr\"odinger systems, especially why solitons are so common in certain regimes. (They are so highly probable that they are the only structures observed macroscopically.)

Consider a wavefunction $u$, that is,  $u$ mapping a time variable $t\in \R$ and space variable $x\in \R^d$ into the complex plane $\C$, which satisfies a $d$-dimensional nonlinear Schr\"odinger equation (NLS): 
\begin{equation}\label{NLS}
i \partial_t u = - \Delta u + \kappa |u|^{p-1}u,
\end{equation}
where $\Delta$ is the Laplacian operator in $\R^d$, $p$ is the nonlinearity parameter, and $\kappa$ is a coupling parameter which is either $+1$ (repulsive/defocusing interaction) or $-1$ (attractive/focusing). As a matter of notation, we will usually write a continuous wavefunction in the form $u(t,x)$ and a discrete wavefunction in the form $f_k(t)= f(t,k)$ where $k$ is the discrete spatial variable. 

We will be interested in chiefly the cubic ($p=3$) and focusing ($\kappa = -1$) NLS, and especially its discrete counterpart via the finite difference approximation, where the spatial variable $k$ is in a discrete set. Also of interest is a modification of \eqref{NLS} with a fractional power of the Laplacian as a non-local dispersive operator, an equation arising naturally in biophysics models of macro-molecules like DNA~\cite{GMCR1, GMCR2, MCGJR}. 

The plan of this article is largely chronological, and background material for the continuous and discrete NLS can be found in~\cite{APT,C,SS}. In section \ref{dNLS}, we will discuss solitons and long-time results for the discrete NLS based on analogues of dispersive PDE methods. Then in section \ref{fNLS} there are similar but more nascent results for discrete NLS-type systems with long-range interactions, again adapted from PDE methods. We introduce the Gibbs measures for continuous NLS in section \ref{SM} and discuss applications to questions of well-posedness for low-regularity initial data (\ref{SMWP}), applications to the typicality of solitons in the absence of wave collapse (\ref{SMS}), and a controversial conjecture for the Gibbs measures of the NLS (\ref{SMconj}). There is a way to get around the problem of nonexistence of Gibbs measures in the presence of wave collapse, so in section \ref{PT} we discuss new results for the Gibbs measures of the 3D discrete NLS, where there is a striking phase transition to soliton-like behavior.

\section{Solitons and breathers for the discrete NLS}\label{dNLS}

Consider a lattice $V$ which is $\Z^d$, for the present, with lattice spacing $h$, i.e., a positive real number denoting the ``distance'' between two neighboring vertices in $V$ \cite{W}. 

The discrete NLS on $V$ is a family of coupled~ODEs governing the evolution of discrete wavefunctions $f_k = f_k(t)$, i.e., maps from $k \in V$ and $t \in \R$ into $\C$ satisfying:
\begin{equation}\label{DNLS}
 i \ddt  f_k  = - \widetilde{\Delta} f_k + \kappa \abs{f_k}^{p-1} f_k, \ \ k\in V.
\end{equation}
Here the discrete nearest-neighbor Laplacian $\widetilde{\Delta}$ is defined by
\[
\widetilde{\Delta} f_k := \frac{1}{h^2}\sum_{j\sim k} (f_j-f_k).
\]  
Note that the sum is over all vertices $j$ that are neighbors of $k$, and that division by $h^2$ means that the discrete Laplacian converges to the continuous Laplacian as $h$ goes to zero. 

Discrete NLS equations can exhibit {\it localized modes} or {\it discrete breathers}, which are periodic-in-time and spatially localized solutions. Breathers have a long history in the literature, and Flach, Kladko, and MacKay made a physics conjecture that the energies of breathers must be bounded away from zero in supercritical lattice models (e.g., discrete NLS in 3D) \cite{FKM}. Weinstein proved one of the conjectures in the supercritical case, that ground state standing waves exist if and only if the total power is larger than some strictly positive threshold~\cite{W}, called the excitation threshold, which is characterized by variational methods.

In order to make this precise, we define the {\it total power} or  {\it particle number} of a wavefunction $f$ solving the discrete NLS \eqref{DNLS} as: 
\[N(f) := || f ||^2_{l^2} = \sum_{k\in V} |f_k|^2.\] 
Also the Hamiltonian energy is well-known to be conserved in time:
\[ H(f) :=  - ( \widetilde{\Delta} f , f) - \frac{2}{p+1} \sum_{k \in V} \abs{f_k}^{p+1} . \]
Here the inner product $(\cdot, \cdot)$ is the usual one on $l^2(\Z^d)$. 
The ansatz for breathers, or standing waves, is that they are of the form for $\omega \in \R$:
\[ f_k(t) = e^{-i \omega t} g_x, \quad k\in V, \quad f_k(t) \in l^2(V).\]
Then $g$ solves the discrete Euler-Lagrange equation
\begin{equation}\label{EL}
 \omega g_x =  - \widetilde{\Delta} g_x - \abs{g_x}^{p-1} g_x.
 \end{equation}
And the solution of \eqref{EL}, called the {\it ground state}, can be constructed as a minimizer of a variational problem:
\[I_{\nu} = \inf \{ H(f):N(f) = \nu \}. \]
Then for $I_{\nu} < 0$, the minimum is attained by $g$ (which is not necessarily unique), with some $\omega = \omega(\nu)$ such that both the constraint ${N}(g) = \nu$ and the Euler-Lagrange system \eqref{EL} hold. Then the main result of \cite{W} is the following.

\begin{thm}[Weinstein \cite{W}]
\label{W}
For supercritical nonlinearities $p\ge \frac4d + 1$ (for example, cubic nonlinearity in three dimensions), there exists an $l^2$-excitation threshold $\nu_c = \nu_c(d, p)>0$ for the ground state of the discrete NLS, meaning that: 
\[I_\nu < 0 \text{ if and only if } \nu> \nu_c.\]

On the other hand, for subcritical $1<p< \frac4d + 1$ (cubic in one dimension), there is no excitation threshold, i.e., $I_\nu<0$ for all $\nu>0$.
\end{thm}

As a preview of the statistical mechanics perspective we will take later, Weinstein's result means that there exist localized modes at zero temperature, corresponding to the ground state. (In the context of statistical mechanics, temperature is a scaling parameter that determines how likely states of various energies are.) Using statistical mechanics, we will be able to add that such localized modes exist for a whole range of temperatures and moreover are typical in a probabilistic sense. Moreover, there will be a phase transition as the mass and the temperature are varied: For low values of mass or inverse temperature, the states that are highly probable are small states with negligible energies. But then there's a sudden jump at a critical threshold of mass or inverse temperature to highly probable localized states with negative energies, in fact bounded away from zero.

Related work of Malomed and Weinstein used this same variational approach to analyze the dynamics of solitons in the 1D discrete NLS: For nonlinearities up to about cubic, there is a unique soliton of large width and no excitation threshold. For nonlinearities in a medium range between cubic and quintic, there is no threshold and there are three solitons, one unstable and two stable (in sharp contrast to the situation for the continuous NLS). And for nonlinearities above quintic there is an excitation threshold for the solitons, one of which is stable and highly localized, and the other unstable and wide \cite{MW}.

Weinstein also made a conjecture in \cite{W} about the subthreshold behavior, that if $\mathcal{N}(f) < \nu_c$ for the solution $f$ of the discrete NLS \eqref{DNLS}, then $f$ disperses to zero, in the sense that: 
\[||f(t)||_{l^q(V)} \xrightarrow{|t| \to \infty} 0,  \text{ for any } q \in (2, \infty].\]
This conjecture has been addressed by similar adaptations of PDE methods to the discrete NLS by Stefanov and Kevrekidis, who proved the following decay and Strichartz-type estimates for discrete Schr\"odinger (and Klein-Gordon) equations.

\begin{thm}[Stefanov-Kevrekidis \cite{SK}]
\label{SK}
For the free discrete Schr\"odinger equation $  i \ddt  f_k  = - \widetilde{\Delta} f_k $ with initial condition $f(0) = f_k(0)$, we have the decay estimate for the wavefunction $f(t) = f_k(t)$:
\[ || f(t) ||_{l^\infty} \le C (1 + |t|)^{-d/3} || f(0) ||_{l^1}.\]
And for the nonlinear discrete Schr\"odinger equation \eqref{DNLS}, with initial data $f(0)$ in the sequence space $l^2$, we have the Strichartz-type estimates for any admissible pairs $(q,r)$ and $(\tilde{q}, \tilde{r})$ (admissible means $2 \le(q,r) \le \infty$,  $\frac1q + \frac{d}{3r} \le \frac{d}{6}$, and $(q,r,d) \ne (2, \infty, 3)$), that
\[ ||f(t)||_{L^q l^r} \le C( ||f(0)||_{l^2} + ||\, |f|^2 f  ||_{L^{\tilde{q}'} l^{\tilde{r}'}}),\]
where $q'$ is the dual exponent of $q$, and the mixed Lebesgue spaces $L^q_t l^r_x$ are defined by:
\[ || f_k(t) ||_{L^q l^r} := \left( \int_0^\infty \left(\sum_{k\in V} |f_k(t)|^r \right)^{q/r} dt \right)^{1/q}. \]
Moreover, the constant in the bound can be chosen uniformly in time for any compact interval $[0,T]$.
\end{thm}

The decay here is not as fast as for the continuous Schr\"odinger decay estimates (which is of order $t^{-d/2}$), but it is sharp in the discrete setting. One consequence of this approach is the following result about decay of small solutions, addressing Weinstein's conjecture and showing that standing wave solutions are impossible in the subthreshold regime.

\begin{corollary}[Stefanov-Kevrekidis \cite{SK}]
\label{SKcor}
If $p\ge \frac6d+1$, then there exists a positive $\eps$ depending on dimension only, such that for small enough initial data, $||f(0)||_{l^2} \le \eps$, a unique solution exists for all time and vanishes in $l^r$ in the limit as $t \to \infty$, for all $r>2$. In fact:
\[ ||f||_{L^q l^r} \le C \eps. \]
Moreover, under certain conditions ($p> \frac4d+1$ and $d \le 2$), small solutions decay like the solution to the free Schr\"odinger equation above: for all $q \in [2,(8-2d)/(d+1)],$\[ ||f(t)||_{l^q} \le C t^{-d(q-2)/(3q)} || f(0)||_{l^{q'}}.\]
\end{corollary}
This is a scattering result, where the solution of the nonlinear equation is asymptotic to ({\it scatters to}) the solution of the linear (free) Schr\"odinger equation as time goes to infinity \cite{SS}. But this uniform smallness of the solution for small initial data happens only when the nonlinearity is more than just supercritical ($\frac6d+1$ rather than $\frac4d+1$).

There are two scaling limits of the discrete Schr\"odinger equation that are of interest. One is the {\it continuum limit} alluded to earlier, of sending the lattice spacing $h$ to zero. In this limit, at least in one and two dimensions, one obtains the continuous nonlinear Schr\"odinger equation \eqref{NLS}. Another limit is the {\it anti-integrable} or {\it anti-continuum} limit, resulting in an infinite system of decoupled ODEs by sending $h$ to infinity, for instance, or equivalently sending the total power $\nu$ to infinity. In this case, ground states become concentrated on the lattice at a single site \cite{W}, and their spectral properties and stability issues have been studied extensively \cite{K3,KPS, PK, PS, PS2}. 

(An interesting thing will happen for the statistical mechanics of the continuum limit in three dimensions, however, because the nonlinearity is supercritical there and thus dominates the kinetic energy, making the coupling effectively negligible and producing behavior similar to the anti-continuum limit in that states typically localize at a single site.)

For instance, Kevrekidis, Pelinovsky, and Stefanov \cite{KPS} build on the previous dispersive decay estimates of \cite{SK}, in order to prove asymptotic stability  of the ground states under some assumptions on the spectrum of the linearization of the DNLS equation with external potential in 1D:
\begin{equation}\label{potDNLS}
 i \ddt  f_k  = (- \widetilde{\Delta} +V_k) f_k + \kappa \abs{f_k}^{p-1} f_k, \ \ k\in V,
\end{equation}
Because of the slower rate of decay in the $l^\infty$ norm of solutions to the discrete Schr\"odinger equations \ref{SK}, the critical power of the nonlinearity is $p= 7$ in place of the continuous Schr\"odinger equation's critical power $p=5$, and \cite{KPS} does not handle the cubic case $p=3$. 

Define the spaces $l^1_s$ and $l^2_s$ by the norms:
\[ ||f||_{l^1_s} := \sum_{k \in \Z} (1+k^2)^{s/2} |f_k|,\]
\[ ||f||_{l^2_s} := \left( \sum_{k \in \Z} (1+k^2)^{s} |f_k|^2 \right)^{1/2}.\]
Decompose the solution into a family of solitons and a radiation part:
\[ u(t) = e^{-i\theta(t)}[ \phi(\omega(t)) + z(t)].\]
\begin{thm}[Kevrekidis-Pelinovsky-Stefanov \cite{KPS}]\label{KPS}
Assume that $\kappa = 1$, that $p \ge 4$, that $V \in l^1_{2\sigma}$ for some $\sigma > \frac52$, that no solution $\psi_0$ of $H\psi_0 = 0$ exists in $l^2_{-s}$ for $s \in (\frac12, \frac32]$, and that $H: = -\tilde{\Delta} + V$ has exactly one negative simple eigenvalue $\omega_0 < 0$ with eigenvector $\psi_0 \in l^2$ and no eigenvalues greater than $4$. Also assume that there is no initial phase shift $\theta(0) = 0$, that the initial soliton location $\omega(0)$ is close to a fixed $\omega_0$, and that the initial data $u(0)$ is close to the soliton $\phi(\omega(0))$ in the space $l^2$. 

Then in the limit as $t$ goes to infinity, $\omega(t)$ converges, and the radiation part $z(t)$ vanishes in $l^\infty$.
\end{thm}

This is a nice asymptotic stability result, but what's more remarkable is that something similar happens in a continuum limit of high-dimensional ($d \ge 3$)  finite lattices called the {\it thermodynamic} limit, sending $h$ to zero and simultaneously taking the total mass (power divided by number of lattice sites) to be fixed (i.e., sending power to infinity in particular way). In the thermodynamic limit, one can find a critical threshold above which one obtains states concentrated states at a single site, at least with high probability.

\section{Solitons for discrete NLS with long-range interactions}\label{fNLS}
In pure physics applications like optics and Bose-Einstein condensation, the NLS is the model, but in the context of biophysics where particle interactions can be long-range, the natural model turns out to be a Schr\"odinger equation with fractional-order Laplacian. An example is modeling electron transport in DNA \cite{GMCR1, GMCR2, MCGJR}. On DNA, and other macromolecules with complicated helical structures, electrons or excitations can propagate along the strand in a very non-local way. 

Mathematically, these phenomena are just beginning to be explored. The usual questions about well-posedness and solitons are interesting for the fractional NLS, and in fact, there is a corresponding discrete system whose solitons have been studied numerically and whose continuum limit is the predicted fractional NLS. Work on the fractional PDEs has been done by
\cite{CSS, GW}, on fractional NLS specifically by \cite{BL,FQT}, and on the related Benjamin-Ono equation by \cite{KMR, M}.

The discrete set-up is on a 1D lattice $h \Z$ with mesh size $h > 0$, and some fixed $h_0$ such that $h < h_0 \leq 1$. Writing $x_k = hk$ with $k \in \Z$, we consider a family of discrete wavefunctions $f^h_k = f^h(t,x_k) : \R \times h \Z \to \C$ that satisfy the following discrete NLS-type equation:
\begin{equation}\label{discrete}
i \frac{d}{dt} f^h_k(t) = h \sum_{j: j \neq k} \frac{ f^h_k(t) - f^h_j(t) }{|x_k-x_j|^{1+2 s} }  +\kappa |f^h_k(t)|^{p-1} f^h_k(t) 
\end{equation}
Here the sum is over all vertices $j$ distinct from fixed $k$, and $0 < s < \infty$ is a fixed parameter governing the decay of the long-range interaction. These wavefunctions have conserved quantities similar to the traditional discrete NLS for energy:
\[
H(f^h) = \frac{h^2}{4} \sum_{j,k: j \neq k} \frac{|f^h_k - f^h_j|^{2}}{|x_k -x_j|^{1+2s}} + \frac{\kappa h}{p+1} \sum_{k} |f^h_k|^{p+1},
\]
and for mass:
\[
M(f^h) =  h \sum_{k} |f^h_k|^2.
\]
Equivalent to this last, power is conserved:
\[
N(f^h) =  \sum_{k} |f^h_k|^2.
\]

As a model of DNA, the cubic nonlinearity describes a self-interaction for a base pair of the strand with itself, and the summation term models interactions between base pairs decaying like an inverse power of the distance along the strand \cite{MCGJR}. The complex twisted structure of DNA in three dimensions is what allows distant base pairs to interact, even at a long distance.

Under reasonable conditions, solutions of this ``fractional" discrete Schr\"odinger equation $\eqref{discrete}$ exist for all time and are unique by standard fixed-point methods \cite{BL,FQT,KLS}, and solitary wave solutions of this system of highly coupled ODEs have been derived \cite{GMCR1} and extended to more general nonlinearities \cite{BL}. 

Explicit solitary wave solutions for this discrete fractional NLS \eqref{discrete} have been studied in the biophysics literature. For instance, Gaididei et al. \cite{GMCR1} make the ansatz $f^h(x_k) = e^{i \lambda t} \phi^h(x_k) $, with localized state:
\begin{equation}\label{fdSoliton}
\phi^h(x_k) = \sqrt{ \frac{h^{-1} N \sinh a}{\cosh (a(2 \delta - 1))}} \exp{(-a |k-\delta|)},
\end{equation}
the parameter $a$ representing the inverse ``width'' of the soliton and $\delta \in [0,1)$ determining its location. Numerically studying the associated variational problem, and minimizing $H$ under the power constraint $N(f^h) = constant$, they find two kinds of localized states, one called ``on-site" with $\delta = 0$ and one ``inter-site" with $\delta = 1/2$, depending on the value of the inverse power parameter $s$. 

For on-site states with cubic nonlinearity, there appears to be a critical value $s_c$, estimated to be near $1$, above which the dispersive interaction decays fast enough that the behavior is qualitatively like that of the nearest-neighbor approximation, or traditional discrete NLS. And below the critical value $s_c$, they find three bistable localized states, one low-frequency wide state and two high-frequency narrow states  \cite{GMCR1}. On the other hand, inter-site states with cubic nonlinearity are not so well understood, and they find a lower critical value $s'_c \simeq .55$ and a smaller regime $s \in (2, s'_c)$ in which there are two stable inter-site states. In both cases, the dependence of the power on the frequency of the localized states is studied.

Gaididei et al. also present the continuum version of these nonlocal models, with $u=u(t,x)$ solving the fractional NLS of the form:
\begin{equation}\label{continuum}
i \partial_t u = c(-\Delta)^{\alpha} u \pm |u|^{p-1} u
\end{equation}
with $u : \R \times \R \to \C$, a constant $c$ depending only on $s$, and $\alpha$ depending on $s$ appropriately. In the biophysics literature, this limit $h \to 0$ is routinely taken without mathematically rigorous arguments, but rigor is forthcoming with the result as follows \cite{KLS}.


We consider the cubic interaction $p = 3$ for simplicity, but this can be generalized to other nonlinearities $p$, and we consider the inverse power law decay in  \eqref{discrete}, but this can be generalized to other interaction kernels in \cite{KLS} that fall into asymptotic classes with behavior substantively similar to the kernel in with parameter $s$. 

The result is that if we are given a family of initial data for the discrete evolution problems \eqref{discrete}, $0<h<h_0$, that approximate a sufficiently smooth initial datum for the fractional NLS \eqref{continuum}, then the solutions of the discrete equations converge in an appropriate weak sense to the solution of the fractional NLS, with the fractional power of the Laplacian:
 \[ \alpha = \alpha(s) = \left \{ \begin{array}{ll} s, & \quad \mbox{for $\frac 1 2 < s < 1$}, \\ 1 , & \quad \mbox{for $s \geq 1$} . \end{array} \right . \]

In other words, the critical value is exactly $s_c = 1$ (not too far from the numerical results in \cite{GMCR1}) below which the long-range interactions produce a nonlocal fractional NLS in the continuum limit, with Laplacian of order $\alpha = s$. Above the critical threshold, the interaction strength decays so quickly that only local effects survive in the continuum limit, which is exactly the ``classical" local NLS, $\alpha = 1$. And at the threshold, we get the classical NLS in the continuum limit, with a logarithmic factor appearing in some scaling constants. 

For this continuum limit result, we develop fractional discrete Sobolev-type inequalities in the discrete setting, uniform embedding, interpolation, and {\it a priori} estimates, and finally for the limit itself, suitable weak compactness arguments. It would be nice also to obtain rigorous results about the solitary wave states from the biophysics literature and their stability for both the discrete and continuous long-range systems, and these long-range NLS systems should be amenable to the statistical mechanics approach discussed next.

\section{Statistical mechanics and Gibbs measures for the continuous NLS}\label{SM}

Adapting tools from the theory of PDEs like variational methods and Strichartz estimates has been useful for understanding the discrete NLS and its localized soliton states, but it is nice to get beyond standard existence results to generic behavior. By setting up the Gibbs measure, or statistical ensemble, which is a probability measure on the space of all solutions, we can tell which behaviors are typical in a probabilistic sense.

{\it Gibbs measures}, or {\it invariant measures}, are based on normalizing the Boltzmann distribution, so that the probability of finding the system in state $u$ is: 
\[\frac{1}{Z(\beta)} e^{-\beta H(u)}.\] 
Here $\beta$ is the inverse temperature, sometimes replaced by $(k_B T)^{-1}$, which makes temperature $T$ explicit using the Boltzmann constant  $k_B$. And $Z(\beta)$ is the normalizing constant, also known as the partition function, which encodes statistical properties of the system (in particular the free energy as we'll see in section \ref{ZF}). For dynamical systems such as the NLS, the energy $H$ is conserved by the dynamics, so that we can invoke the Liouville theorem to see that the measure is invariant.

The intuition behind the Boltzmann distribution can be understood by thinking about the extreme cases: In the high temperature case, $\beta \ll 1$, the Gibbs measure is essentially uniform, that is, the system could statistically be in almost any state. On the other hand, near zero temperature, $\beta \gg 1$, the Gibbs measure is nearly vanishing for any state but the lowest-energy one, and thus the system has a strong preference for the ground state(s). The statistical ensemble associated to the Gibbs measure is called the {\it canonical ensemble}, and the states that the system is likely to be in, according to the Gibbs measure, are called {\it canonical macrostates}, so called because those are the states that are observed macroscopically.

There are variations on the standard Gibbs measure and its canonical ensemble. For instance, the {\it microcanonical ensemble} arises from conditioning the Gibbs measure on the energy $H$ being a constant $E$ (really being within $\eps$ of $E$ and sending $\eps$ to zero appropriately).
A natural question to ask is whether the microcanonical ensemble is equivalent to the canonical ensemble, meaning that they share the same macrostates. The answer is yes and no: if certain conditions are satisfied, they are equivalent, and if not, the microcanonical ensemble has a richer set of macrostates. Another ensemble is the {\it petit canonical} ensemble, which results from conditioning on mass: this is also called a {\it mixed ensemble} because it is canonical in energy and microcanonical in mass.

An advantage of Gibbs measures is that they ignore transient states and focus instead on the long-time probabilistic behavior of the system, and this approach was initiated in the work of Lebowitz, Rose, and Speer \cite{LRS} for the 1D focusing NLS on the circle $\T$:
\begin{equation}\label{1DNLS}
 i \partial_t u = - \partial_{xx} u - \abs{u}^{p-1} u. 
\end{equation}
 The canonical ensemble isn't the right thing to look at because it is not normalizable, so they conjectured that, with a mass cutoff, one can construct a formal invariant measure from the Hamiltonian energy:
\begin{equation}\label{cont}
H(u) : = \frac{1}{2} \int_0^L \abs{u ' (x)}^2 dx - \frac{1}{p+1} \int_0^L \abs{u  (x)}^{p+1} dx.
\end{equation}
The reason for the mass cutoff is that if one forms the Gibbs measure with density $e^{-\beta H(u)}$, then it is not normalizable because $Z(\beta)$ is infinite. This can be fixed by making $B$ the allowed mass: the partition function $Z$ is then a (finite) function of both the thermodynamic $\beta$ and the allowed mass $B$. Thus we form the following formal invariant measure: 
\[ (Z(\beta,B))^{-1} e^{-\beta H(u)}  \mathbf{1}_{\{\lVert u \rVert^2_2 \leq B \}} \prod_x du(x).\]
Here the product needs to be interpreted properly, highlighting one of the difficulties in constructing Gibbs measures for infinite-dimensional systems: a finite-dimensional approximation must be done, for instance, by spectral truncation or by spatial discretization.

\subsection{Application of Gibbs measures to well-posedness}\label{SMWP}
Bourgain continued the investigation by using spectral truncation to avoid the harmonic analysis on cyclic groups that would have been necessary for spatial discretization \cite{B1}. He introduced a discretization in Fourier space by truncating the Fourier sum at the $n$-th Fourier coefficient:
\[ u^n := P_n u := \sum_{|k| \leq n} \hat{u}(k) e^{2\pi i k x}.\]
Then the finite-dimensional Hamiltonian system is:
\begin{equation}\label{stNLS} i \partial_t u^n = - \partial_{xx} u^n - P_n\left(\abs{u^n}^{p-1} u^n \right),\end{equation} 
with the Hamiltonian energy functional:
\[ H(u^n) = 2 \pi^2 \sum_{|k| \leq n} k^2 |a_k|^2 - \frac{1}{p+1} \int_{\T} \Bigg| \sum_{|k| \leq n} a_k e^{2 \pi i kx} \Bigg|^{p+1} dx,\]
where $a_k$ is the $k$-th Fourier coefficient: $ a_k := \hat{u}(k).$
The Hamiltonian is conserved, as is the truncated mass: 
\[\left( \sum_{|k| \leq n} |a_k|^2 \right)^{\frac12}.\]
Next consider a ball in this finite-dimensional phase space:
\[\Omega_{n,B} : = \left\{ (a_k)_{|k|\leq n} \Bigg| \sum_{|k| \leq n} |a_k|^2  \leq B^2 \right\},\]
and the natural Gibbs measure on $\Omega_{n,B}$ at fixed inverse temperature $\beta = 1$ is then:
\[ d\mu_n :=  \exp{\left\{ \frac{1}{p+1}  \int_{\T} \Bigg| \sum_{|k| \leq n} a_k e^{2 \pi i kx} \Bigg|^{p+1} dx\right\}} \rho_n \otimes da_0,\]
where $\rho_n := Z_n^{-1} \exp{\{- 2 \pi^2 \sum_{0< |k| \leq n} k^2 |a_k|^2\}} $ and the normalization constant,
\[Z_n := \int_{\C} \cdots \int_{\C} e^{- 2 \pi^2 \sum k^2 |a_k|^2} da_1 da_{-1} \cdots da_n da_{-n}.\] 
We note that this is a petit, or mixed, canonical ensemble: canonical in $H$ and microcanonical in the mass, i.e., the $L^2$ norm is constrained to be no more than a constant $B$. The limit of $\rho_n$ as $n \to \infty$ is $\rho$, which turns out to be the image measure of the random Fourier series map:
\[ \omega \mapsto \sum_{k \neq 0} \frac{g_k(\omega)}{2 \pi k} e^{2 \pi i kx}.\]
The $g_k$ are independent identically distributed complex Gaussian variables, and the image is almost surely in $H^{s}$ for any $s<\frac12$. And the limit of $\mu_n$ as $n \to \infty$ is $\mu$, the following weighted Wiener measure:
\begin{equation}\label{mu}
d\mu = Z^{-1} e^{\frac{1}{4} \lVert u \rVert^4_4} \mathbf{1}_{\{\lVert u \rVert^2_2 \leq B \}} d a_o \otimes d \rho.
\end{equation}
Here $d\rho$ can be represented either through the random Fourier series leaving off the zero mode, or the formal Weiner measure (with its own normalization $Z_0$):
\[ d a_o \otimes d \rho = Z_0^{-1} e^{- \frac12 \int_\T |\partial_x u|^2 dx} \prod_{x \in \T} du(x).\]

Then the invariance of the measures $\mu_n$ under the finite-dimensional flow $S_n$ of the spectrally truncated NLS \eqref{stNLS}, combined with the uniform regularity of $S_n$ in $H^s$ and the approximation of the NLS flow by $S_n$, all give the following theorem for the limit measure: 
\begin{thm}[Bourgain \cite{B1}]
The measure $\mu$ in \eqref{mu} is invariant under the dynamics of the periodic 1D focusing cubic ($p = 3$) and quintic ($p = 5$) NLS.
\end{thm}
This is also parlayed into a global well-posedness result for the continuous NLS:
\begin{thm}[Bourgain \cite{B1}]
The Cauchy problem for the 1D periodic NLS with $p \leq 5$ and with random initial data:
\begin{equation*}\begin{cases} i \partial_t u = - \partial_{xx} u - \abs{u}^{p-1} u \\
u(0,x) = a_0 + \sum_{k \neq 0} \frac{g_k(\omega)}{2 \pi k} e^{2 \pi i kx},
\end{cases} \end{equation*}
is globally well-posed for almost all $a_0$ and $\omega$ satisfying the mass constraint:
\[ |a_0|^2 + \sum_{k \neq 0} \frac{|g_k(\omega)|^2}{4 \pi^2 k^2} \leq B^2.\]
\end{thm}

Another advantage of Gibbs measures has been illustrated in recent work on the 1D cubic NLS with low-regularity initial data on the circle $\T$, namely that in spite of the ill-posedness results for the NLS with initial data of lower regularity than $L^2$ \cite{CCT}, the bad initial data are quite rare, forming a set of vanishing measure with respect to the natural Gibbs measure. In fact, for the cubic NLS on the circle, Colliander and Oh were able to prove almost-sure local well-posedness for initial data in $H^\sigma(\T)$, for every $\sigma>-1/3$, and global well-posedness for data in $H^s(\T)$, for $s > -1/12$ \cite{CO}. Oh and Sulem proved weak continuity of the cubic NLS solution map from $L^2(\T)$ to the space of distributions \cite{OS}. Both papers relied on a Wick ordering of the nonlinearity, which seems to be an appropriate approach for the 1D cubic NLS with rough initial data; by contrast Molinet proved that the standard solution map (without Wick ordering) from $L^2(\T)$ to the space of distributions is not weakly continuous. This approach has also been fruitful for almost sure global well-posedness for the periodic derivative NLS \cite{NORBS}. 

\subsection{Application of Gibbs measures to solitons}\label{SMS}
Jordan, Josserand, and Turkington have established the connection between Gibbs measures and solitons for the 1D NLS with non-traditional (i.e., not polynomial) nonlinearities \cite{JJ,JT}. 
Consider the following NLS on a finite interval $[0,L]$ with periodic boundary conditions:
\begin{equation}\label{JTNLS}
 i \partial_t u = - \partial_{xx} u - g(|u|^2) u,
\end{equation}
with focusing and bounded nonlinearity, i.e., $g(0) = 0, \; \forall a \ge 0, g'(a)> 0$ and $g(a)\leq C$, and as $a \to \infty, ag'(a) \leq C'$; for example a nonlinearity that approximates behavior of large wave amplitudes:
\[g(|u|^2) = \frac{|u|^2}{1 + |u|^2}.\]
As usual the two conserved quantities are energy:
\[ H(u) := \frac12 \int_0^L (|\partial_x u|^2 - G(|u|^2)) dx,\]
where $G$ is the primitive of $g$: $G(a) = \int_0^a g(s) ds$, and mass:
\[ M(u) := \frac12 \int_0^L |u|^2 dx.\]

In this 1D setting, the energy appearing in the Gibbs measure contains the kinetic term which captures the fine-grained fluctuations and the potential term which essentially ignores those fluctuations.


Using the same kind of spectral truncation as Bourgain, the mixed ensemble measure is formed with a regularization with respect to mass as well as the usual conditioning on constrained mass:
\[P^B_\beta(du^n):= (Z(\beta_n,B))^{-1} \exp{(\beta[H(u^n) + \sigma \Vert u^n\Vert^2_{2} ])}  \mathbf{1}_{\{\lVert u ^n\rVert^2_2 \leq B \}} \prod du^n.\]
They fix the mass $B$ and scale the inverse temperature as $\beta_n \sim Cn$ for large $n$ so that the mean energy has a finite limit, then rigorously analyze the limit using large deviations theory, comparing the result with numerical simulations \cite{JJ}. The large deviations principle is that for appropriate sets $S \subset L^2$:
\[\lim_{n \to \infty} \frac{1}{\beta_n} \log P^B_{\beta_n} \{u^{n} \in B \} = - \inf_{u \in S} I^B(u),\]
where $I^B$ is the large deviations rate function defined for $u \in H^1_0$:
\[ I^B(u) := 
\begin{cases}
H(u) - E(N), \;  \frac12 \Vert u\Vert^2_{2} = B,\\
\infty, \; \text{ otherwise.}
\end{cases}\]
Here $E$ is the coherent energy function $E(N):= \min \{ H(u) : M(u) = B, u \in H^1_0 \}.$

What this means is that with high probability the random fields (spectrally truncated wavefunctions) $u^{n}$ converge exactly to the ground state solitons of the continuous NLS. This is because the minimizers of the large deviations rate function are the states that are most probable macroscopically, and all other states have exponentially vanishing probability. And since the large deviations rate function is equal to the Hamiltonian energy up to a constant, the ground state solitons that solve the traditional variational problem are identical with the most probable macroscopic states. 

The ground states also turn out to be stable both macroscopically and microscopically: the large deviations principle says that the macrostate stays close to the ground states with high probability, even if some of the energy goes into turbulent fluctuations; and the typical microstate as sampled from the Gibbs measure $P^B_{\beta_n}$ is close to a ground state in the $H^1$-norm.

Two key things that Jordan, Josserand, and Turkington need in order for their approach to work are non-integrability and the absence of wave collapse, so in particular, it doesn't work for the cubic NLS in one, two, or three dimensions. It turns out that the 1D and 3D cases can be handled in the right frameworks, and possibly two dimensions as well.

\subsection{Conjecture about a phase transition to soliton-like behavior}\label{SMconj}
Numerical simulations led Lebowitz, Rose, and Speer to make a second conjecture about a possible phase transition in the invariant measures for the 1D NLS: at high temperature and low mass they saw uniformly small wavefunctions, but at low temperature and high mass there emerged sharply concentrated structures like solitons.


By contrast, simulations in \cite{Burl} suggested that there is no phase transition, and then Rider, following on the works of McKean and Vaninsky \cite{mckean94, mckean97a, mckean97b}, confirmed this for the 1D infinite-volume focusing NLS by proving that the thermodynamic limit is trivial \cite{Rider1, Rider2}.  Bourgain also studied invariant measures of the 1D infinite-volume defocusing NLS \cite{B2} and of the 2D defocusing NLS \cite{B3} (see also the review article \cite{B4}), as did Tzvetkov \cite{T}. 

These statistical mechanics methods that work so well for the subcritical NLS break down at criticality. For instance, Brydges and Slade studied the 2D focusing NLS \cite{BS} and saw that the natural Gibbs measure construction cannot produce an invariant measure for large coupling coefficients of the (critical) nonlinearity. Similarly for 3D focusing NLS, the natural construction is not normalizable, and it is thought to be impossible to make any reasonable Gibbs construction.

This is the motivation for the alternative approach explained below, using a spatial discretization to construct the finite-dimensional approximation corresponding to the focusing cubic NLS in three dimensions. There are still unresolved issues about what information can be transferred to the continuum limit, but this approach has turned out to be very fruitful, especially as far as how much can be said about properties of the Gibbs measures--as it turns out, a lot more than the previous spectral approaches, and even about a phase transition to solitons.

\section{A phase transition to solitons for 3D discrete NLS}\label{PT}

Because constructing Gibbs measures for the 3D NLS appears to be impossible, we consider instead the NLS on a spatial discretization of the 3D unit torus $[0,1]^3$, represented by: 
\[V = \{0, 1/L,\ldots, (L-1)/L\}^3,\]
so that each vertex of $V$ has degree $d=3$, $V$ has $n:= |V| = L^3$ sites, and the inter-site distance is $h= 1/L$.

(We consider this interesting special case here for simplicity, but in \cite{CK} the results actually handle general graphs $G = (V,E)$ that are finite and undirected without self-loops, and that are {\it high-dimensional} in the sense that $nh^2$ goes to $\infty$ as $h \to 0$, where $h>0$ is the distance between two neighboring vertices in $G$ and $n$ is the number of sites, $n \to \infty$. This condition holds for the discrete 3D torus defined above, but does not hold for the discrete 2D torus in the usual continuum limit.)

The focusing~NLS on the discrete 3D torus is:
\begin{equation}\label{DNLSlam}
 i \frac{d}{dt}  f_k  = - \widetilde{\Delta} f_k - \abs{f_k}^2 f_k,
 \end{equation}
with discrete Laplacian:
\begin{equation}\label{DL} \widetilde{\Delta} f_k = \frac{1}{h^2}\sum_{j: j\sim k} (f_j-f_k) .\end{equation}

Solutions of \eqref{DNLSlam} are known to exist for all time \cite{W}, and the dynamics of this system conserves the power $N(f) := \sum_{k\in V} |f_k|^2$, or equivalently mass $M(f):=n^{-1} N(f)$, and the energy, appropriately normalized by the number of sites $n$:
\[ H(f) := \frac{2}{n} \sum_{\substack{ j,k \in V \\ k \sim j}} \biggl|\frac{f_k - f_j}{h}\biggr|^2 - \frac{1}{n}\sum_{k\in V} |f_k|^4. \]
Then by the Liouville theorem, the natural Gibbs measure $d\mu := e^{- \beta H(f) } \prod_{k\in V} d f_k$ is invariant under the dynamics of the discrete NLS \eqref{DNLSlam} for any real $\beta$ (inverse temperature). However, this measure has infinite mass if $\beta > 0$, so we construct the Gibbs measures with a mass cut-off to restrict the system to an allowed mass $B$:
\begin{equation}\label{tmu}
d\tilde{\mu}_{\beta,B} := Z^{-1} e^{- \beta H(f) }  1_{\{N(f) \le Bn\}} \prod_{k\in V} d f_k.
\end{equation}
Here $B$ is arbitrary and positive, and the normalizing constant, or partition function, is $Z = Z(\beta,B)$. This cut-off Gibbs measure is still invariant under the NLS dynamics because mass is a conserved quantity, and it turns out to be tractable in the limit as the grid size $h$ goes to zero, at least in dimensions three and higher. 

Let $\psi$ be a random element of $\C^V$ with probability distribution $\tilde{\mu}$: that is, $\psi$ is a random wavefunction on $V$ such that for each $A \subseteq \C^V$,
\[
\pp(\psi \in A) = Z^{-1} \int_{A} e^{-\beta H(f) } 1_{\{N(f) \le Bn\}} df, 
\]
where $df = \prod_x df_k$ denotes the Lebesgue differential element on $\C^V$. In order to understand the behavior of the random map $\psi$, we first study the partition function $Z$. The first theorem below shows that if we have a sequence of graphs with $n$ and $nh^2$ both tending to infinity, the limit of  $n^{-1}\log Z$ can be exactly computed for any positive $\beta$ and $B$. 

\subsection{The partition function and free energy}\label{ZF}

\begin{figure}[t!]
      \centering
      \includegraphics[height=2.6in,clip, angle=270]{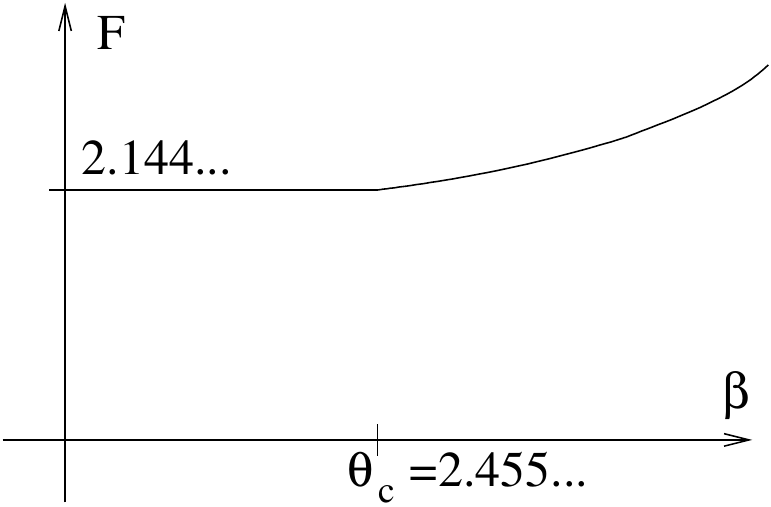} 
       \caption{The free energy is constant for small inverse temperature and starts increasing at the critical threshold. Here mass is normalized, $B=1$.}
       \label{fig1}
\end{figure}

The first result is that for positive temperature, and for the discrete torus in dimensions three and higher, there exists a function $F$ that we construct, such that as $n \to \infty$:
\[ Z(\beta,B) \sim e^{nF(\beta,B)}.\] 
More precisely, let $g:[2,\infty) \to \R$ be the function
\begin{equation}\label{mdef}
g(\theta) := \frac{\theta}{2} - \frac{1}{2}+ \frac{\theta}{2} \sqrt{1-\frac{2}{\theta}} + \log\biggl(\frac{1}{2}-\frac{1}{2}\sqrt{1-\frac{2}{\theta}}\biggr). 
\end{equation}
It can be checked that $g$ has a unique real zero that we call $\theta_c$, because $g$ is strictly increasing in $[2,\infty)$, $g(2) <0$ and $g(3) >0$.  Numerically, $\theta_c \approx 2.455407$. Let\begin{equation}\label{fbeta}
F(\beta, B) := 
\begin{cases}
\log(B\pi e) & \text{ if } \beta B^2 \le \theta_c,\\
\log(B\pi e) + g(\beta B^2) &\text{ if } \beta B^2 > \theta_c.
\end{cases} 
\end{equation}
(Figure \ref{fig1} shows  a graph of $F$ versus $\beta$ when $B=1$.)
The theorem also gives an explicit rate of convergence. 
\begin{thm}[Chatterjee-Kirkpatrick \cite{CK}]
\label{zthm}
Suppose $\beta \ge 0$. Let $\ep\in (0,1/5)$. There exists a positive constant $C$ depending only on $\ep$, $\beta$, $B$, $h$ and $d$ such that if $n > C$, then 
\[
\frac{\log Z}{n} - F(\beta, B) \ge -Cn^{-1/5}- C(nh^2)^{-1}
\]
and 
\[
\frac{\log Z}{n} - F(\beta, B)\le 
\begin{cases}
Cn^{-1/5+\ep} + C n^{-4\ep/5}  &\text{ if } \beta B^2 \le \theta_c,\\ 
Cn^{-1/5+\ep}  &\text{ if } \beta B^2 > \theta_c.
\end{cases}
\]
\end{thm}
The proof of these asymptotics for the partition function is at the heart of the work \cite{CK}: large deviations techniques for discrete random variables are used, first obtaining upper and lower bounds for $\frac{\log Z}{n}$ and then optimizing which gives the free energy $F$. Then, because $F$ is not differentiable at the critical threshold $\theta_c$, that means there is a first-order (discontinuous) phase transition at  $\theta_c$. The random map $\psi$ behaves quite differently in the two phases, and there is an important quantity which is discontinuous at the threshold.

\subsection{The phase transition and its consequences}

\begin{figure}[htpb]
      \centering
        \includegraphics[height=2.2in,clip, angle=270]{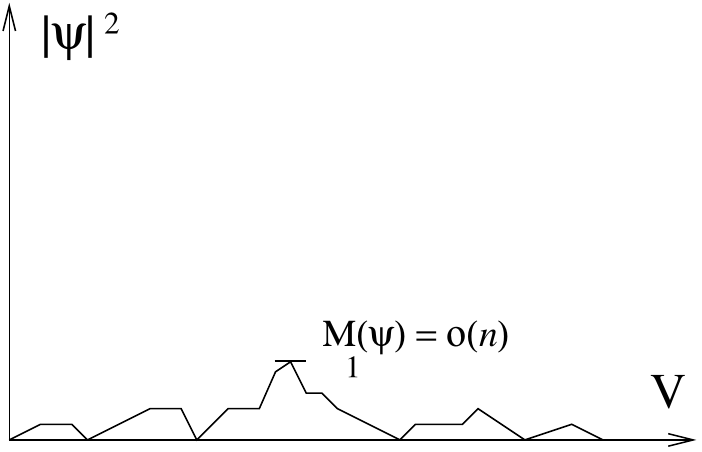}   
      \includegraphics[height=2.2in,clip, angle=270]{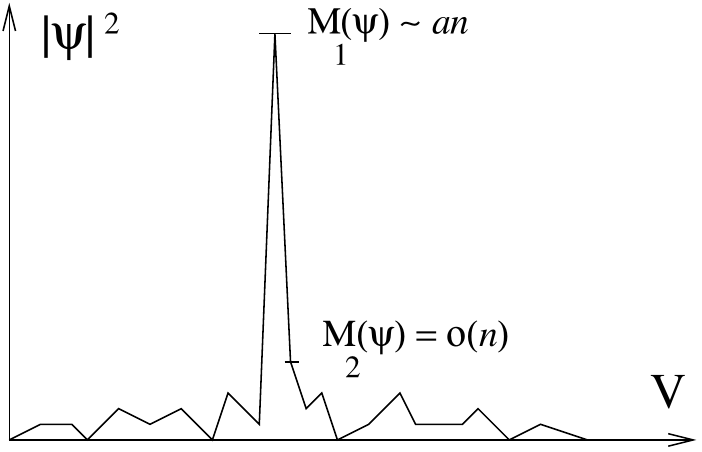}
      \caption{Impressionistic pictures of typical wavefunctions in the two phases. Left: subcritical noncondensed phase, $\beta B^2< \theta_c$, where the heaviest site bears a negligible proportion of the mass compared to the total mass $N(\psi) \sim Bn$. Right: supercritical condensed phase, $\beta B^2> \theta_c$, where a significant proportion $a$ of the mass concentrates at a single site.}
      \label{impression}
\end{figure}

To roughly describe this phase transition (depicted in figure \ref{impression}), when $\beta B^2> \theta_c$ there is a single site $k\in V$ which bears an abnormally large fraction of the total mass of the random wavefunction $\psi$. This fraction is nearly deterministic, given by the ratio $a/B$, where 
\begin{equation}\label{adef}
a = a(\beta, B) := \frac{B}{2}+\frac{B}{2}\sqrt{1-\frac{2}{\beta B^2}}.
\end{equation}

More precisely, if $M_1(\psi)$ and $M_2(\psi)$ are the largest and second largest components of the vector $(|\psi_k|^2)_{k\in V}$, then with high probability $M_1(\psi)\approx an$ and $M_2(\psi) = o(n)$. Moreover, $N(\psi) \approx Bn$ with high probability.  A consequence is that the largest component carries more than half of the total mass: 
\[
\max_k\frac{|\psi_k|^2}{\sum_j |\psi_j|^2} \approx \frac{a}{B} >\frac{1}{2}. 
\]
On the other hand, when $\beta B^2 < \theta_c$, then $M_1(\psi) = o(n)$, but still $N(\psi)\approx Bn$. Consequently
\[
\max_k\frac{|\psi_k|^2}{\sum_j |\psi_j|^2} \approx 0.
\]
\begin{figure}[t!]
      \centering
      \includegraphics[height=2.5in,clip, angle=270]{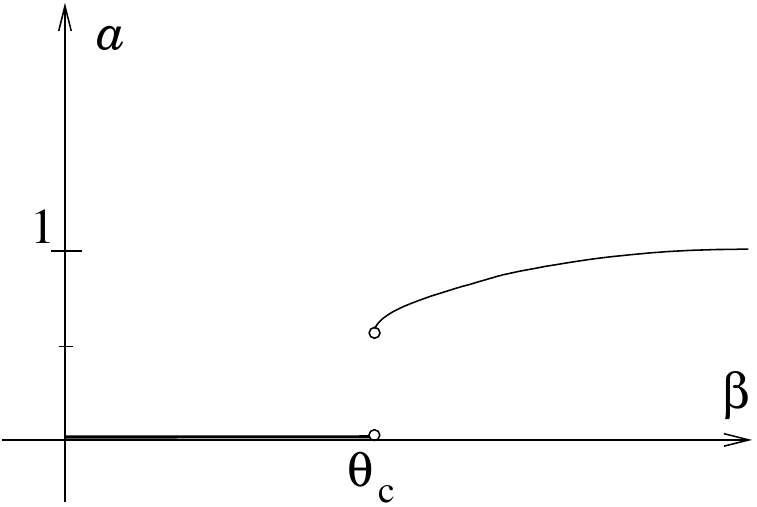} 
       \caption{The fraction of mass at the heaviest site jumps from roughly zero for small inverse temperature, to roughly $.71$ at the critical threshold. (Here $B=1$.)}
       \label{fig2}
\end{figure}

The formula for $a$ (see also the graph of $a$ in figure \ref{fig2} as a function of $\beta$ for $B=1$) shows that $a$ does not tend to zero as $\beta B^2$ approaches $\theta_c$ from above, making it clear that this is a discontinuous (in physics parlance, first-order) phase transition.

Also, in the regime $\beta B^2 < \theta_c$, the energy density is close to zero. By contrast, when $\beta B^2 > \theta_c$, the energy density $H(\psi)/n$ is strictly negative and approximately equals $-a^2$, the exceptional site bearing nearly {\it all} of the energy of the system. This is because the total energy $H(\psi)$ is approximately $-a^2 n$, while the energy at $k$ is, summing over just the neighbors $j$ of $k$: 
\begin{equation*}\begin{split}
\frac{1}{nh^2}\sum_{j\sim k} |f_k-f_j|^2 - \frac{|f_k|^4}{n} & \approx -n^{-1}M_1(\psi)^2 + O(h^{-2}) \\
& = -a^2 n + o(n),
\end{split}
\end{equation*}
the equality by Theorem \ref{phase}, which makes precise this phase transition:

\begin{thm}[Chatterjee-Kirkpatrick \cite{CK}]
\label{phase}
Let $V$ be the discrete torus in dimension $d \geq 3$, let $a=a(\beta, B)$ be defined as in \eqref{adef}, and let $M_1(\psi)$ and $M_2(\psi)$ be the largest and second largest components of $(|\psi_k|^2)_{k \in V}$. 
\begin{enumerate}
\item First, suppose $\beta B^2 >\theta_c$. 
For any $q \in (\frac45, 1)$, there is a constant $C$ depending only on $\beta$, $B$, $d$, and $q$  such that if $n >C$, then with probability $\ge 1- e^{-n^q/C}$, 
\begin{equation}\label{phase1}
\begin{split}
&\biggl|\frac{H(\psi)}{n}+ a^2\biggr| \le Cn^{-(1-q)/4}, \ \ \ \biggl|\frac{N(\psi)}{n} - B\biggr| \le Cn^{-(1-q)/2}, \\
&\biggl|\frac{M_1(\psi)}{n} - a\biggr|\le Cn^{-(1-q)/4}, \ \text{ and} \ \ \frac{M_2(\psi)}{n}\le  Cn^{-(1-q)}. 
\end{split}
\end{equation}
\item Next, suppose $\beta B^2 < \theta_c$. 
For any $q \in (\frac{17}{18},1)  $, there is a constant $C$ depending only on $\beta$, $B$, $d$, and $q$ such that whenever $n>C$, with probability $\ge 1- e^{-n^q/C}$ 
\begin{equation}\label{phase2}
\begin{split}
&\biggl|\frac{H(\psi)}{n}\biggr| \le 2n^{-2(1-q)}, \ \ 
\biggl|\frac{N(\psi)}{n} - B\biggr| \le n^{-(1-q)}, \\&\text{and } \ \frac{M_1(\psi)}{n} \le n^{-(1-q)}.
\end{split}
\end{equation}
\item Finally, if $\beta B^2=\theta_c$ and $q \in (\frac{17}{18},1) $,
then there is a constant $C$ depending only on $\beta$, $B$, $d$, $p$ and $q$ such that whenever $n>C$, with probability $\ge 1- e^{-n^q/C}$ either \eqref{phase1} or \eqref{phase2} holds. 
\end{enumerate}
\end{thm}
The proof proceeds by a multi-step approximation: first the Gibbs measures are approximated by uniform measures on appropriate spaces (intersections of $l^2$ and $l^4$ spheres), and then the uniform measures are approximated by i.i.d. complex Gaussian random variables on $l^4$ annuli. Then the i.i.d. Gaussians must concentrate on the parts of the $l^4$ annuli closest to the origin, i.e., where one coordinate is large and the others are uniformly small.

\begin{figure}[htpb]
      \centering       \includegraphics[height=2.5in,clip, angle=270]{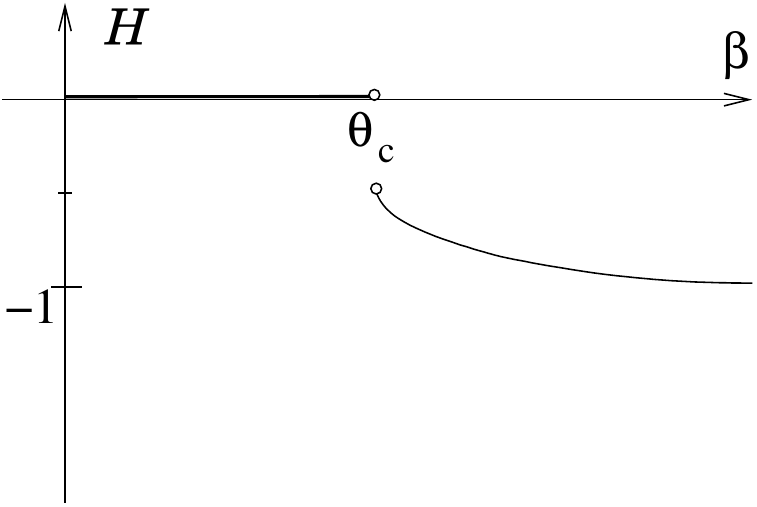} 
       \caption{Localized states occur only when the energy is nonzero, and in fact bounded away from zero, as in the physics conjecture.}
       \label{fig4}
\end{figure}

One consequence of this theorem is confirmation of the conjecture of Flach, Kladko, and MacKay, because the energies of the localized states are bounded away from zero, made clear in figure \ref{fig4} where energy is plotted as a function of $\beta$.

Another straightforward consequence of this theorem is that typical discrete wavefunctions above the critical threshold blow up in the discrete $H^1$ norm: 
$$\|f\|^2_{\widetilde{H}^1} :=\frac{1}{n}\sum_{k\in V} |f_k|^2 +  \frac{1}{n} \sum_{(k,j) \in E} \biggl| \frac{ f_k - f_j }{h}\biggr|^2.$$
It also turns out that the discrete $H^1$ norm also diverges even for $\beta B^2 \le \theta_c$, just with different rates:
\begin{thm}[Chatterjee-Kirkpatrick \cite{CK}]\label{blowup}
Suppose the context of Theorem \ref{phase} holds, in dimension $d \geq 3$, with $n > C$. Let $\psi = (\psi_k)_{k \in V}$ be a discrete wavefunction picked randomly from the invariant probability measure $\tilde{\mu}$ defined in~\eqref{tmu}. If $\beta B^2 >\theta_c$, then there is a positive constant $c$ depending only on $\beta$, $B$, $d$ and $p$ such that $\pp(\| \psi \|_{\widetilde{H}^1} \le cn^{p})\le e^{-n^c}$ whenever $n \ge 1/c$. On the other hand, if $\beta B^2 \le \theta_c$, then similarly: $\pp(\| \psi \|_{\widetilde{H}^1}\le c\sqrt{d} n^{p})\le e^{-d n^c}$. 
\end{thm}

Since the measure $\tilde{\mu}$ of \eqref{tmu} is invariant for the discrete NLS \eqref{DNLSlam}, one may expect that if the initial data comes from $\tilde{\mu}$, localized modes will continue to exist for all time. The question is whether the mode jumps around or stays in one place (in which case we have a standing or stationary wave, also called a discrete breather). The following theorem uses the statistical equilibrium result to deduce a dynamical result: the same site continues to be the localized mode for an exponentially long period of time. 
\begin{thm}[Chatterjee-Kirkpatrick \cite{CK}]\label{expo}
Suppose the previous context holds with $\beta B^2 > \theta_c$ and the function $a$ as defined in \eqref{adef}. Let $\psi(t) = (\psi_k(t))_{k\in V}$ be a discrete wavefunction evolving according to \eqref{DNLSlam}, where the initial data $\psi(0)$ is picked randomly from the invariant probability measure $\tilde{\mu}$ defined in~\eqref{tmu}. Then for any $q \in (\frac45,1)$,
there is a constant $C$ depending only on $\beta$, $B$, $d$, and $q$ such that if $n > C$, then with probability $\ge 1-e^{-n^q/C}$ the inequalities \eqref{phase1} hold for $\psi(t)$ for all $0 \le t \le e^{n^q/C}$, and moreover there is a single $k\in V$ such that the maximum of $|\psi_j(t)|$ is attained at $y=x$ for all~$0\le t\le e^{n^q/C}$. In particular, $\psi(t)$ is approximately  a standing wave with localized mode at $k$ for an exponentially long time. 
\end{thm}
This dynamical stability result combines PDE and probabilistic methods. First, with high probability the wavefunction is condensed at time $t=0$, and then the invariance of the Gibbs measure implies that the wavefunction is condensed for an overwhelming proportion of times up to the exponential time in the theorem. If the soliton were to jump from one site to another, then it would have to happen in a tiny window of time; but that would result in a wavefunction that violates basic NLS estimates, so in fact the soliton must stay put.

In other words, it is not only possible, but typical, for solutions of \eqref{DNLSlam} to have unique stable localized modes for exponentially long times if the initial inverse temperature or mass are above a threshold. Earlier results for localized modes were only for existence, and only for large mass, i.e.,\ tending to infinity (e.g.,\ \cite{W}), while Theorem \ref{expo} proves typicality for finite mass.

Our final result is about the probability  distribution of the individual coordinates of a random map $\psi$ picked from the measure $\tilde{\mu}_{\beta, B}$. We can give a precise description of the distribution for small collections of coordinates: In both phases, the joint distributions are approximately complex Gaussian vectors with the appropriate variance ($B$ below the threshold and $B-a$ above it, with $a$ defined in \eqref{adef}), provided the collection is sufficiently small compared to the number of sites in the grid and provided a symmetry assumption is satisfied.

The symmetry assumption on $G$ is that $G$ is {\it translatable} by some group~$\Sigma$: where $\Sigma$ is a group of automorphisms of $G$ such that:  1) $|\Sigma|=n$, and 2) no element of $\Sigma$ except the identity has any fixed point. For example, the discrete torus is translatable by the group of translations.

\begin{thm}[Chatterjee-Kirkpatrick \cite{CK}]\label{distribution}
Suppose the graph $G$ is translatable by some group of automorphisms, $d \geq 3$, and $\psi$ is a random wavefunction picked according to the measure $\tilde{\mu}$. Take any $m$ distinct points $x_1,\ldots, x_m\in V$. Let $\phi = (\phi_1,\ldots, \phi_m)$ be a vector of i.i.d.\ standard complex Gaussian random variables. If $\beta B^2 < \theta_c$, then there is a constant $C>0$ depending only on $\beta$, $B$, and $d$ such that if $n > C$, then for all Borel sets $U\subseteq \C^m$,
\[
\bigl|\pp(B^{-1/2}(\psi_{x_1},\ldots,\psi_{x_m}) \in U) - \pp(\phi\in U)\bigr|\le mn^{-1/C}. 
\]
If $\beta B^2 > \theta_c$, the result holds after $B^{-1/2}$ is replaced with $(B-a)^{-1/2}$ where $a = a(\beta, B)$ is defined in \eqref{adef}, and the error bound is changed to $m^3 n^{-1/C}$. 
\end{thm}
This says that below the critical threshold, the wavefunctions look like white noise, but above the threshold they look like a single large coordinate and white noise elsewhere.

As for the future, one direction to push these kinds of results is into the two-dimensional setting. There the nonlinear and kinetic energies are in a delicate balance, and the wavefunction can't be approximated by i.i.d. Gaussians; instead one must consider correlations between sites, producing what is called a discrete Gaussian free field. Another direction is to study the Gibbs measures for the fractional NLS, which will have non-local behavior coming from the underlying L\'evy processes. There is also a possibility that this statistical mechanics approach could be useful for the soliton resolution conjecture: that in certain regimes, solutions generically tend asymptotically towards a soliton plus the small fluctuations of a radiation component.

\bibliographystyle{amsalpha} 

\thispagestyle{plain}

\end{document}